\documentclass[conference,a4paper]{APSIPA2018}
\usepackage{multirow}
\usepackage{graphicx}
\usepackage{amsmath}
\usepackage[psamsfonts]{amssymb}
\usepackage{amsxtra}
\usepackage{makecell}
\usepackage{romannum}
\usepackage{threeparttable}
\usepackage{lipsum}

\newcommand{\tensor}[1]{\boldsymbol{\mathcal{#1}}}

\newcommand{\mat}[1]{\mathbf{#1}}
\newcommand{\vect}[1]{\mathbf{#1}}

\makeatletter

\newcommand{\Rmnum}[1]{\expandafter\@slowromancap\romannumeral #1@}
\makeatother
\renewcommand{\figurename}{Figure} 
\renewcommand{\tablename}{Table} 
\renewcommand{\thefootnote}{\alph{footnote}}
\begin{document}
\title{Higher-dimension Tensor Completion via Low-rank Tensor Ring Decomposition}

\author{%
\authorblockN{%
Longhao Yuan\authorrefmark{1}\authorrefmark{2}, Jianting Cao\authorrefmark{4}\authorrefmark{1}, Xuyang Zhao\authorrefmark{1}\authorrefmark{2}, Qiang Wu\authorrefmark{5}\authorrefmark{6} and Qibin Zhao\authorrefmark{2}\authorrefmark{3}
}
\authorblockA{%
\authorrefmark{1}
Graduate School of Engineering, Saitama Institute of Technology, Japan}
\authorblockA{%
\authorrefmark{2}
Tensor Learning Unit, RIKEN Center for Advanced Intelligence Project (AIP), Japan}
\authorblockA{%
\authorrefmark{3}
School of Automation, Guangdong University of Technology, China}
\authorblockA{%
\authorrefmark{4}
School of Computer Science and Technology, Hangzhou Dianzi University, China}
\authorblockA{%
\authorrefmark{5}
School of Information Science and Engineering, Shandong University, China}
\authorblockA{%
\authorrefmark{6}
Institute of Brain and Brain-Inspired Science, Shandong University, China}
}

\maketitle
\newcommand\blfootnote[1]{%
\begingroup 
\renewcommand\thefootnote{}\footnote{#1}%
\addtocounter{footnote}{-1}%
\endgroup 
}
\thispagestyle{empty}

\begin{abstract}
The problem of incomplete data is common in signal processing and machine learning. Tensor completion algorithms aim to recover the incomplete data from its partially observed entries. In this paper, taking advantages of high compressibility and flexibility of recently proposed tensor ring (TR) decomposition, we propose a new tensor completion approach named tensor ring weighted optimization (TR-WOPT). It finds the latent factors of the incomplete tensor by gradient descent algorithm, then the latent factors are employed to predict the missing entries of the tensor. We conduct various tensor completion experiments on synthetic data and real-world data. The simulation results show that TR-WOPT performs well in various high-dimension tensors. Furthermore, image completion results show that our proposed algorithm outperforms the state-of-the-art algorithms in many situations. Especially when the missing rate of the test images is high (e.g., over 0.9), the performance of our TR-WOPT is significantly better than the compared algorithms.
\end{abstract}
\blfootnote{Corresponding authors: Jianting Cao (cao@sit.ac.jp) and Qibin Zhao (qibin.zhao@riken.jp).}
\section{Introduction}
Tensors are high-dimension representations of vectors and matrices. Many kinds of data in real-world, for example, color images ($length\times width\times RGB \ channels$), videos ($length\times width \times  RGB \ channels \times time$) and electroencephalography (EEG) signals ($magnitude\times trails \times time$) are more than two dimensions. Usually, the high-dimension data is first transformed to vector or matrix, only then can the data be applied to traditional algorithms. In that way, the adjacent structure information of the original data will be lost, thus lead to redundant space cost and low-efficiency computation \cite{shashua2005non}. Tensor representation can retain the high dimension of data and solve the above problems. Tensor has been studied for more than a century, many methodologies have been proposed \cite{kolda2009tensor}. Moreover, tensor has been applied in various research field such as signal processing \cite{cichocki2015tensor}, machine learning \cite{anandkumar2014tensor}, data completion \cite{acar2011scalable},  brain-computer interface (BCI) \cite{liu2014tensor}, etc..

CANDECOMP/PARAFAC (CP) decomposition \cite{bro1997parafac} and Tucker decomposition \cite{tucker1966some} are the most popular and classic tensor decomposition models, of which many theories and applications have been proposed \cite{sidiropoulos2017tensor}. In recent years, a new theory system named tensor network has drawn people's attention and becomes a promising aspect of tensor methodology \cite{cichocki2016tensor,cichocki2017tensor}. One of the most representative tensor decomposition models of tensor network is matrix product state (MPS), which is also known as tensor train (TT) decomposition \cite{oseledets2011tensor}. TT decomposition shows high data compression ability and computational efficiency. One of the most significant features of TT decomposition is that it can overcome the ``curse of dimensionality", i.e., for an $N$-dimension tensor, the number of parameters of Tucker decomposition is exponential in $N$, while the number of parameters of TT decomposition is linear in $N$. Although CP decomposition achieves high compression ability and the number of parameters is also linear in $N$, finding the optimal latent factors of CP decomposition is very difficult. Recently, a new tensor decomposition model named tensor ring (TR) decomposition which is more generalized than TT decomposition has been proposed \cite{zhao2016tensor}. Tensor ring owns all the dominant properties TT has and it has a more flexible model constraint than tensor train. TR relaxes the rank constraint of TT, thus lead to more interesting properties, e.g., enhanced compression ability, interpretability of latent factors and rotational invariance ability of latent factors \cite{zhao2017learning}. Several papers have studied and applied TR decomposition in machine learning fields \cite{wang2018wide,cao2017tensorizing}, in which TR shows promising aspects in representation ability and computational efficiency. 

Tensor completion is to recover an incomplete tensor from the partially observed entries of the tensor, which has been applied in various completion problems such as image/video completion \cite{liu2013tensor,zhao2016bayesian}, compressed sensing \cite{gandy2011tensor}, link prediction \cite{lu2011link}, recommendation system \cite{karatzoglou2010multiverse}, etc.. The theoretical key point of matrix completion and tensor completion is the low-rank assumption of data. There exists strong theoretical support and various solutions for solving the low-rank problem of matrices, and the most studied convex relaxation of low-rank matrix is nuclear norm \cite{hu2013fast}. However, determining the rank of a tensor is an NP-hard problem \cite{hillar2013most}. To solve this problem, there are mainly two types of tensor completion methods, ``rank minimization based" approach and ``tensor decomposition based" approach \cite{long2018low}. The first approach formulates the convex surrogate models of low-rank tensors. The most representative model of this kind of tensor completion approach employs low-rank constraints on the matricization of every mode of the incomplete tensor \cite{liu2013tensor}. For an $N$-dimension incomplete tensor $\tensor{T}$, the low-rank tensor completion model is formulated by:
\begin{equation}
\label{lr_tensor}
\min \limits_{\tensor{X}} \ \  \sum\nolimits_{n=1}^N\Vert \mat{X}_{(n)} \Vert_* ,\ s.t. \ P_{\Omega}(\tensor{X})=P_{\Omega}(\tensor{T}),
\end{equation}
where $\tensor{X}$ is the low-rank approximation tensor, $\Vert \cdot \Vert_*$ is the nuclear norm, and $P_{\Omega}(\tensor{T})$ denotes the entries w.r.t. the set of indices of observed entries represented by $\Omega$. The missing entries of $\tensor{T}$ is approximated by $\tensor{X}$ and the rank of the completed tensor $\tensor{X}$ is determined automatically. Different from ``rank minimization based" approach, the ``tensor decomposition based" approach do not find the low-rank tensor directly, instead, it firstly finds the tensor decomposition of the incomplete data by observed entries, then the latent factors are used to predict the missing entries. This kind of approach sets the rank of tensor decomposition manually, and the optimization model is given below:
\begin{equation}
\label{lnr_tensor}
\min\limits_{\{\tensor{G}^{(n)}\}_{n=1}^N} \Vert \tensor{W}*(\tensor{T}-\tensor{X}(\{\tensor{G}^{(n)}\}_{n=1}^N) \Vert_F^2,
\end{equation}
where $\Vert \cdot \Vert_F$ is the Frobenius norm, $\{\tensor{G}^{(n)}\}_{n=1}^N$ is the sequence of latent factors under consideration and $\tensor{X}(\{\tensor{G}^{(n)}\}_{n=1}^N)$ is the tensor approximated by the latent factors. $\tensor{W}$ is a weight tensor which is the same size as the incomplete tensor, and every entry of $\tensor{W}$ meets:
\begin{equation}
\label{weight}
 w_{i_{1}i_{2}\cdots i_{N}}=
 \left\{
 \begin{aligned}
 &0 \qquad \text{if} \; y_{i_{1}i_{2}\cdots i_{N}} \;\text{is a missing entry},\\
 &1 \qquad \text{if} \;  y_{i_{1}i_{2}\cdots i_{N}}\;\text{is an observed entry}.
\end{aligned}
\right.
\end{equation}
More explanations of notations can be found in Section $\Rmnum{2}.A$. Based on different tensor decomposition models, various ``tensor decomposition based" approaches have been proposed, e.g., CP weighted optimization \cite{acar2011scalable}, weighted tucker \cite{filipovic2015tucker} and TT weighted optimization \cite{yuan2017completion}. All the methods aim to find the specific structure of the incomplete data by different kinds of tensor decompositions. However, CP, Tucker and TT based WOPT algorithms apply tensor decomposition models which lack of flexibility, this may lead to bad convergence when considering different kinds of data. TR decomposition model is much more flexible, so TR-based WOPT method can get a better approximation of incomplete data, thus provide better completion results.

In this paper, we propose a novel tensor completion method which shows good performance in various experimental situations. The main works of this paper are listed below:
\begin{itemize}
\item Based on the recently proposed TR decomposition, we propose a new tensor completion algorithm named tensor ring weighted optimization (TR-WOPT).

\item The TR latent factors are optimized by gradient descent method and then they are used to predict the missing entries of the incomplete tensor.

\item  We conduct several simulation experiments and real-world data experiments. The experiment results show that our method outperforms the state-of-the-art tensor completion algorithms in various situations. We also find that our method is robust to tensor dimension. When the image data is tensorized to a proper high-dimension, the performance of our method can be enhanced. 
\end{itemize}

\section{Preliminaries}
\subsection{Notations}
In this paper, a scalar is denoted by a normal lowercase or capital letter, e.g., $x \in\mathbb{R}$ and $X \in\mathbb{R}$. A vector is denoted by a boldface lowercase letter, e.g., $\vect{x}\in\mathbb{R}^{I}$. A matrix is denoted by a boldface capital letter, e.g., $\mat{X}\in\mathbb{R}^{I\times J}$. A tensor of dimension $N\geq 3$ is denoted by Euler script letters, e.g., $\tensor{X}\in\mathbb{R}^{I_1\times I_2\times\cdots \times I_N}$. For $n=1,\cdots,N$, $\{ \tensor{X}^{(n)}\}_{n=1}^N:=\{\tensor{X}^{(1)},\tensor{X}^{(2)},\cdots,\tensor{X}^{(N)}\}$ is defined as a tensor sequence, and $\tensor{X}^{(n)}$ is the $n$th tensor of the sequence. The representations of matrix sequence and vector sequence are denoted in the same way. An element of  tensor $\tensor{X}  \in\mathbb{R}^{I_1\times I_2\times\cdots \times I_N}$ of index $\{i_{1},i_{2},\ldots,i_{N}\}$ is denoted by $x_{i_{1}i_{2}\cdots i_{N}}$ or $\tensor{X}(i_{1},i_{2},\ldots,i_{N})$. One type of the mode-$n$ matricization (unfolding) of tensor $\tensor{X}  \in\mathbb{R}^{I_1\times I_2\times\cdots \times I_N}$ is denoted by $\mat{X}_{(n)}\in\mathbb{R}^{I_n \times  {I_1 \cdots I_{n-1} I_{n+1} \cdots I_N}}$. Another type of mode-$n$ matricization of tensor $\tensor{X}  \in\mathbb{R}^{I_1\times I_2\times\cdots \times I_N}$ is denoted by $\mat{X}_{<n>}\in\mathbb{R}^{I_n \times  {I_{n+1} \cdots I_{N} I_{1} \cdots I_{n-1}}}$. The inner product of two tensors $\tensor{X}$, $\tensor{Y}$ with the same size $\mathbb{R}^{I_1\times I_2\times\cdots \times I_N}$ is defined as $\langle \tensor{X},\tensor{Y} \rangle=\sum_{i_1}\sum_{i_2}\cdots\sum_{i_N}x_{i_1 i_2\cdots i_N}y_{i_1 i_2\cdots i_N}$. Furthermore, the Frobenius norm of $\tensor{X}$ is defined by $\left \| \tensor{X} \right \|_F=\sqrt{\langle \tensor{X},\tensor{X} \rangle}$. The Hadamard product is denoted by `$\ast$' and it is an element-wise product of vectors, matrices or tensors of same sizes. For example, given tensor $\tensor{X}, \tensor{Y}\in\mathbb{R}^{I_1\times I_2\times\cdots \times I_N}$, $\tensor{Z}=\tensor{X}*\tensor{Y}$, then $\tensor{Z}\in\mathbb{R}^{I_1\times I_2\times\cdots \times I_N}$, and $z_{i_1 i_2 \cdots i_N}=x_{i_1 i_2 \cdots i_N} y_{i_1 i_2 \cdots i_N}$.

\subsection{Tensor ring decomposition}

\begin{figure}
\centering
\includegraphics[width=0.4\textwidth]{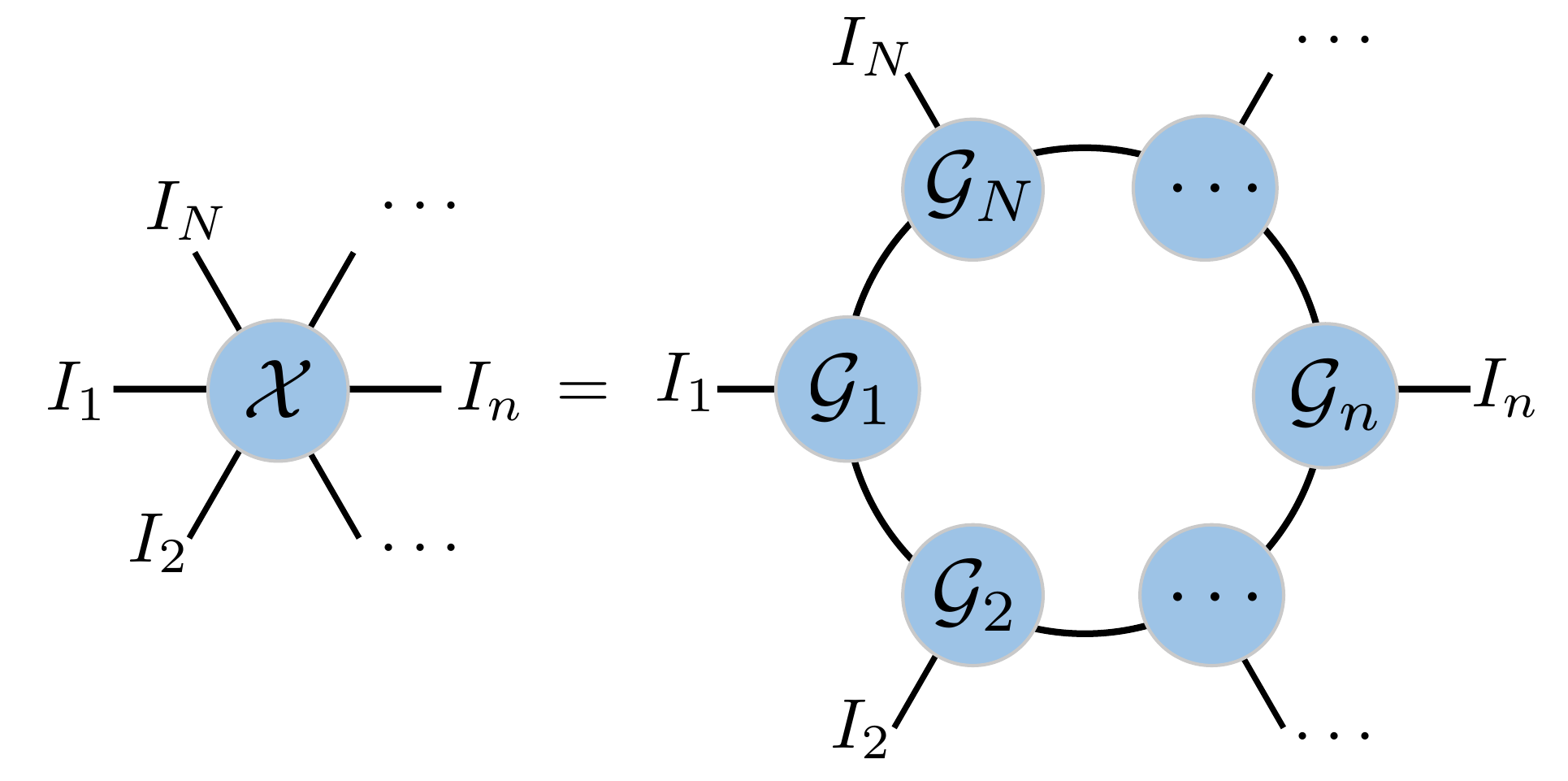}
\caption{TR decomposition.}
\label{tr_d}
\end{figure}

Tensor ring (TR) decomposition is the generalization of tensor-train (TT) decomposition.  It represents a tensor by circular multilinear products over a sequence of latent factors of lower-dimension \cite{zhao2016tensor}. The latent factors are named core tensors. All of the core tensors of TR decomposition are three-dimension tensors, and are denoted by $\tensor{G}^{(n)} \in\mathbb{R}^{R_{n} \times I_{n} \times R_{n+1}}$, $n=1,\ldots,N$, where  $R_1,\ldots,R_N$ denotes the TR-rank which controls the model complexity of TR decomposition. The representation regulation of tensor network in monography \cite{cichocki2016tensor} is adopted and the diagram of TR decomposition is shown in Figure \ref{tr_d}. Similar as TT, the TR decomposition linearly scales the tensor dimension, and thus it can overcome the ``curse of dimensionality''. The TR decomposition relaxes the rank constraints on the first and last core tensors of TT to $R_1=R_{N+1}$, while the original constraints on TT decomposition is rather stringent, i.e., $R_1=R_{N+1}=1$. The rank constraints of TT is because of the relation of core tensors and original tensor has to be deduced by multiple multiplication operation. TR applies trace operation, so all the core tensors can be constrained to be third-dimension equivalently. Therefore, TR is considered as a linear combination of TT and it offers a more powerful and generalized representation ability than TT decomposition. The element-wise relation of the TR core tensors and the original tensor is given by:

\begin{equation}
\label{tr_relation1}
\tensor{X}(i_1,i_2,\ldots,i_N)=\text{Trace}\Big(  \prod_{n=1}^N \mat{G}^{(n)}(i_n) \Big),
\end{equation}
where $\text{Trace}( \cdot )$ is the matrix trace operator, $ \mat{G}^{(n)}(i_n)\in\mathbb{R}^{R_n\times R_{n+1}}$ is the $i_n$th mode-$2$ slice matrix of $\tensor{G}^{(n)}$, which also can be denoted by $\tensor{G}^{(n)}(:,i_n,:)$. Every element of the tensor can be calculated by the trace of the matrices multiple multiplication of the according mode-$2$ slice of every TR core tensors.

\section{Tensor ring weighted optimization}
Based on TR decomposition, we propose the TR-WOPT algorithm which is illustrated as follows. Define $\tensor{T}\in\mathbb{R}^{I_1\times I_2\times\cdots \times I_N}$ is the incomplete tensor with missing entries filled with zero, $\tensor{X}(\{\tensor{G}^{(n)}\}_{n=1}^N)$ is the tensor approximated by the core tensors of TR decomposition. The proposed algorithm is based on ``tensor decomposition based" approach and the model is described in (\ref{lnr_tensor}). The model is to find the core tensors of TR decomposition of an incomplete tensor, then use the TR core tensors to approximate the missing entries. To minimize the model by gradient-based algorithm, the problem is reformulated by the below optimization model:
\begin{equation}
\label{of1}
f(\tensor{G}^{(1)},\ldots,\tensor{G}^{(n)})=\frac{1}{2} \left \|\tensor{W}\ast (\tensor{T}-\tensor{X}(\{\tensor{G}^{(n)}\}_{n=1}^N)) \right \|^{2}_F.
\end{equation}
This is an objective function of an optimization problem and all the core tensors are the optimization objective. From \cite{zhao2016tensor}, the relation between the approximated tensor $\tensor{X}$ and the core tensors $\{\tensor{G}^{(n)}\}_{n=1}^N$ can be deduced as the following equation:
\begin{equation}
\label{XG_relation}
\mat{X}_{<n>}=\mat{G}^{(n)}_{(2)}(\mat{G}_{<2>}^{(\neq n)})^T,
\end{equation}
where $\tensor{G}^{(\neq n)}\in\mathbb{R}^{R_{n+1}\times \prod_{i=1, i\neq n}^N I_i \times R_n}$ is a subchain tensor by merging all core tensors except the $n$th core tensor, i.e., $\tensor{G}^{(\neq n)}:=\tensor{X}(\{\tensor{G}^{(n+1)},\ldots,\tensor{G}^{(n)}, \tensor{G}^{(1)}, \ldots ,\tensor{G}^{(n-1)}\})$. Because each of the core tensors is independent, we can optimize them independently. The optimization function w.r.t. $\tensor{G}^{(n)}$ can be written as:
\begin{equation}
\label{of2}
\begin{aligned}
f(\tensor{G}^{(n)})=\frac{1}{2} \left \|\mat{W}_{<n>}\ast (\mat{T}_{<n>}-\mat{G}_{(2)}^{(n)}(\mat{G}_{<2>}^{(\neq n)})^T) \right \|^{2}_F,
\end{aligned}
\end{equation}
where we consider other tensor cores remain fixed. Next, we can deduce the partial derivatives of the objective function (\ref{of2}) w.r.t. $\mat{G}^{(n)}_{(2)}$ as follow:
\begin{equation}
\label{derTR-WOPT}
\frac{\partial{f}}{\partial{\mat{G}_{(2)}^{(n)}}}=(\mat{W}_{<n>}\ast (\mat{G}_{(2)}^{(n)}(\mat{G}_{<2>}^{(\neq n)})^T-\mat{T}_{<n>})\mat{G}_{<2>}^{(\neq n)}.
\end{equation}
For $n=1,...,N$, the gradients of all the core tensors can be obtained, and the core tensors can be optimized by any gradient-based optimization algorithms. Furthermore, if there is no missing entries in tensor data, our algorithm can also be used as a TR decomposition algorithm. The whole process of applying TR-WOPT to tensor completion is listed in Algorithm 1. The details of gradient descent algorithm and parameter settings in this paper will be explained in Section $\Rmnum{4}$. 

\begin{table}[!htb]
\footnotesize
\begin{center} 
\begin{tabular}{l}
\hline
\textbf{Algorithm 1} Tensor ring weighted optimization (TR-WOPT)\\
\hline
1: \;\;  \textbf{Input}: incomplete tensor $\tensor{T}$, weight tensor $\tensor{W}$, TR-rank \\
 \;\;\;\;\;\;\;\;\;\;\;\;\;\;\;\;\;\;$R_1,\ldots,R_N$, and randomly initialized $\{\tensor{G}^{(n)}\}_{n=1}^N $.\\
2: \;\; \textbf{While} the stopping condition is not satisfied \\
3:  \;\;\;\;\;  \textbf{For} n=1:N \\
4:  \;\; \;\;\;\; Compute gradients of $\{\tensor{G}^{(n)}\}_{n=1}^N $ according to (\ref{derTR-WOPT}).\\
5: \;\;\;\;\; \textbf{End}\\
6: \;\; Update $\{\tensor{G}^{(n)}\}_{n=1}^N $ by gradient descend algorithm. \\
7: \;\; \textbf{End while}\\
8:\;\;\;    $ \tensor{Y}=P_{\Omega}(\tensor{T})+P_{\bar{\Omega}}(\tensor{X}(\{\tensor{G}^{(n)}\}_{n=1}^N ))$\\
9: \;\; \textbf{Output}: completed tensor $\tensor{Y}$.\\
\hline
\end{tabular}
\end{center}
\end{table}

\section{Experiment results}
We conduct several synthetic data experiments and RGB image data experiments to test the performance of our TR-WOPT algorithm. We compare our algorithm with some tensor completion algorithms which are similar to our algorithm, e.g., TT-WOPT \cite{yuan2017completion} and CP-WOPT \cite{acar2011scalable}. Moreover, some other state-of-the-art algorithms, e.g., FBCP \cite{zhao2015bayesian} and FaLRTC \cite{liu2014tensor} are also tested in the next experiments. 

For the gradient descent method applied in our TR-WOPT, we use the nonlinear conjugate gradient (NCG) with line search method, which is implemented by a Matlab toolbox named Poblano toolbox \cite{dunlavy2010poblano}. We use two stopping conditions to all the algorithms: the number of iteration reaches 500 and the error between two iterations satisfies $\Vert\tensor{Y}_2-\tensor{Y}_1\Vert/\Vert\tensor{Y}_2\Vert<tol=10^{-6}$. When one of the stopping conditions is satisfied, the optimization will be stopped.

For completion performance evaluation, we adopt relative square error (RSE) and peak signal-to-noise ratio (PSNR). RSE is calculated by:
\begin{equation}
\text{RSE}=\frac{\Vert  \tensor{T}_{real}-\tensor{Y} \Vert_F}{\Vert \tensor{T}_{real}\Vert_F},
\end{equation}
where $\tensor{T}_{real}$ is the real tensor with full observations, $\tensor{Y}$ is the completed tensor. Moreover, for RGB image data, PSNR is obtained by:
\begin{equation}
\text{PSNR}=10\log_{10}(255^2/\text{MSE}).
\end{equation}
MSE is deduced by:
\begin{equation}
\text{MSE}=\Vert \tensor{T}_{real}-\tensor{Y} \Vert_F^2/\text{num}(\tensor{T}_{real}),
\end{equation}
where num($\cdot$) denotes the number of element of the tensor. In addition, for tensor experiments of random missing cases, we define missing rate as $1-M/\text{num}(\tensor{T}_{real})$, where $M$ is the number of sampled entries (i.e. observed entries).

\subsection{Synthetic data}
\begin{figure}[h]
\centering
\includegraphics[width=0.48\textwidth]{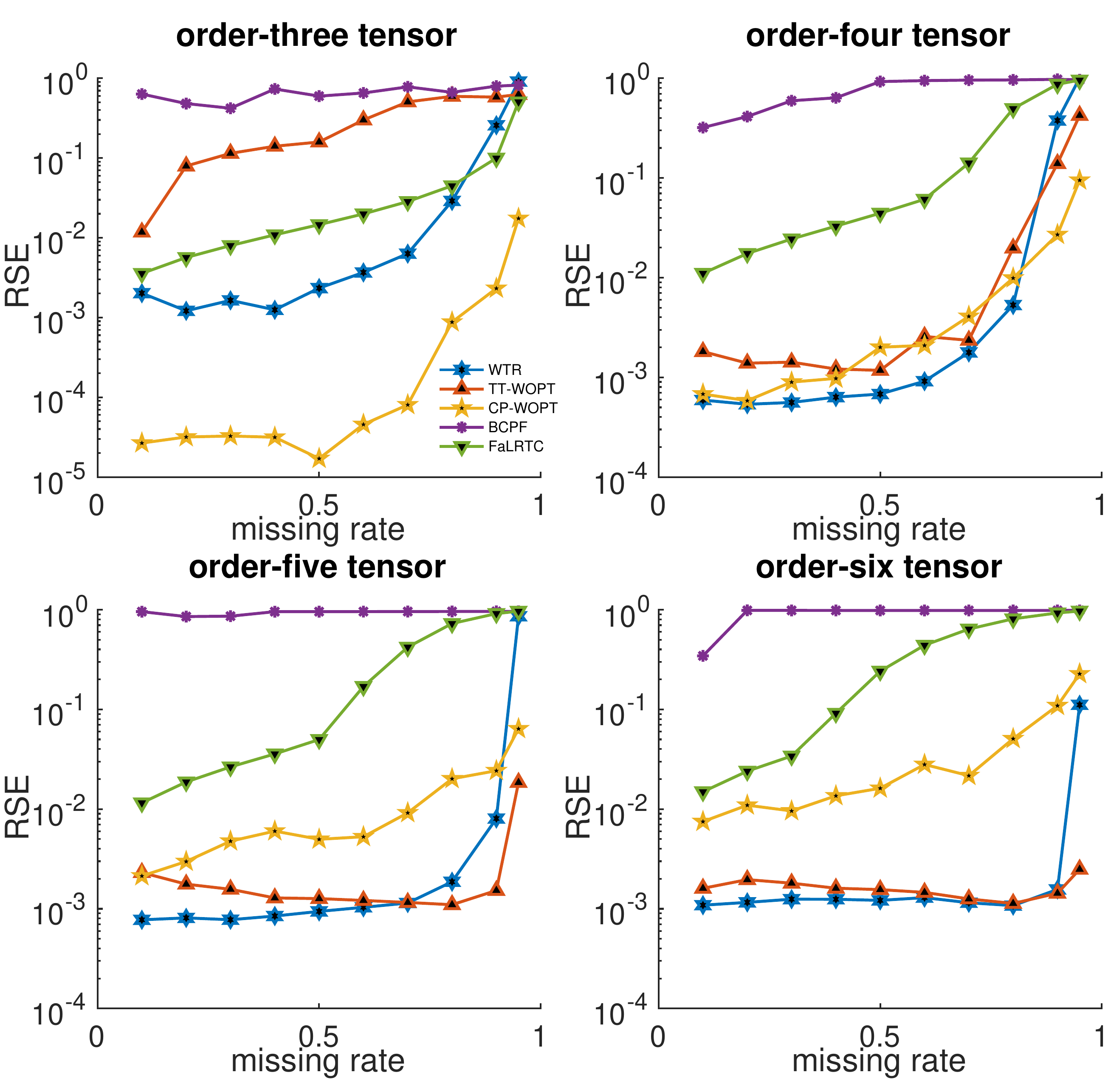}
\caption{Synthetic data completion results of five algorithms under different tensor dimensions and different missing rates.}
\label{sim}
\end{figure}

In this section, we conduct synthetic data experiments to see the performance of our algorithm and the compared algorithms under different tensor dimensions and different missing rates. TR-WOPT, TT-WOPT, CP-WOPT, FBCP, and FaLRTC are tested in the experiments. The synthetic data is generated by a highly oscillating function: $f(x)=sin(\pi/4)cos(x^2)$ in vector form, then the synthetic data is reshaped to generate the required tensors of different dimensions. We test four different tensor dimensions: $48\times48\times48$ (3D), $16\times16\times16\times16$ (4D), $10\times10\times10\times10\times10$ (5D) and $7\times7\times7\times7\times7\times7$ (6D). The missing rate varies from 0.1 to 0.95. For the rank selection of TR-WOPT, TT-WOPT, CP-WOPT, we adjust the ranks of each algorithm to make the model parameters to be as close as possible, thus provide a fair situation to compare model representation ability. FBCP can tune CP rank automatically, and FaLRTC does not need to set tensor ranks beforehand, thus we only need to tune hyper-parameters of the two algorithms and record the best performance. The five algorithms are tested and the RSE results are shown in Figure \ref{sim}. 

From the figure we can see, TR-WOPT performs well in all the four tensor dimensions. Moreover, when the tensor dimension increases, many compared algorithms show performance degradation. TT-WOPT is the closest algorithm to TR-WOPT, though the number of model parameters is set similarly, TR-WOPT performs better than TT-WOPT almost in all situations due to higher representation ability of the TR-based model.

\subsection{Image experiments}

The synthetic data experiments results show that our TR-WOPT performs well in higher-dimension tensors, so in this section, we first test the performance of TR-WOPT under different tensor dimension and TR-ranks. Instead of directly reshaping image to higher-dimension, we apply a novel tensorization scheme mentioned in \cite{yuan2017completion}. The idea is simple, if the size of a RGB image is $U\times V \times 3$, and $U=u_1\times u_2\times \cdots \times u_l$ and $V=v_1\times v_2\times \cdots \times v_l $ are satisfied, then the image can be tensorized to a $(l+1)$-dimension tensor of size $u_1v_1 \times u_2v_2 \times \cdots \times u_lv_l \times 3$. The first dimension of the higher-dimension tensor stands for a pixel block of the image, and the following dimensions are the extension blocks of the image. The higher-dimension tensor generated by the above tensorization scheme is considered to be a better structure of the image.

\begin{figure}
\centering
\includegraphics[width=0.45\textwidth]{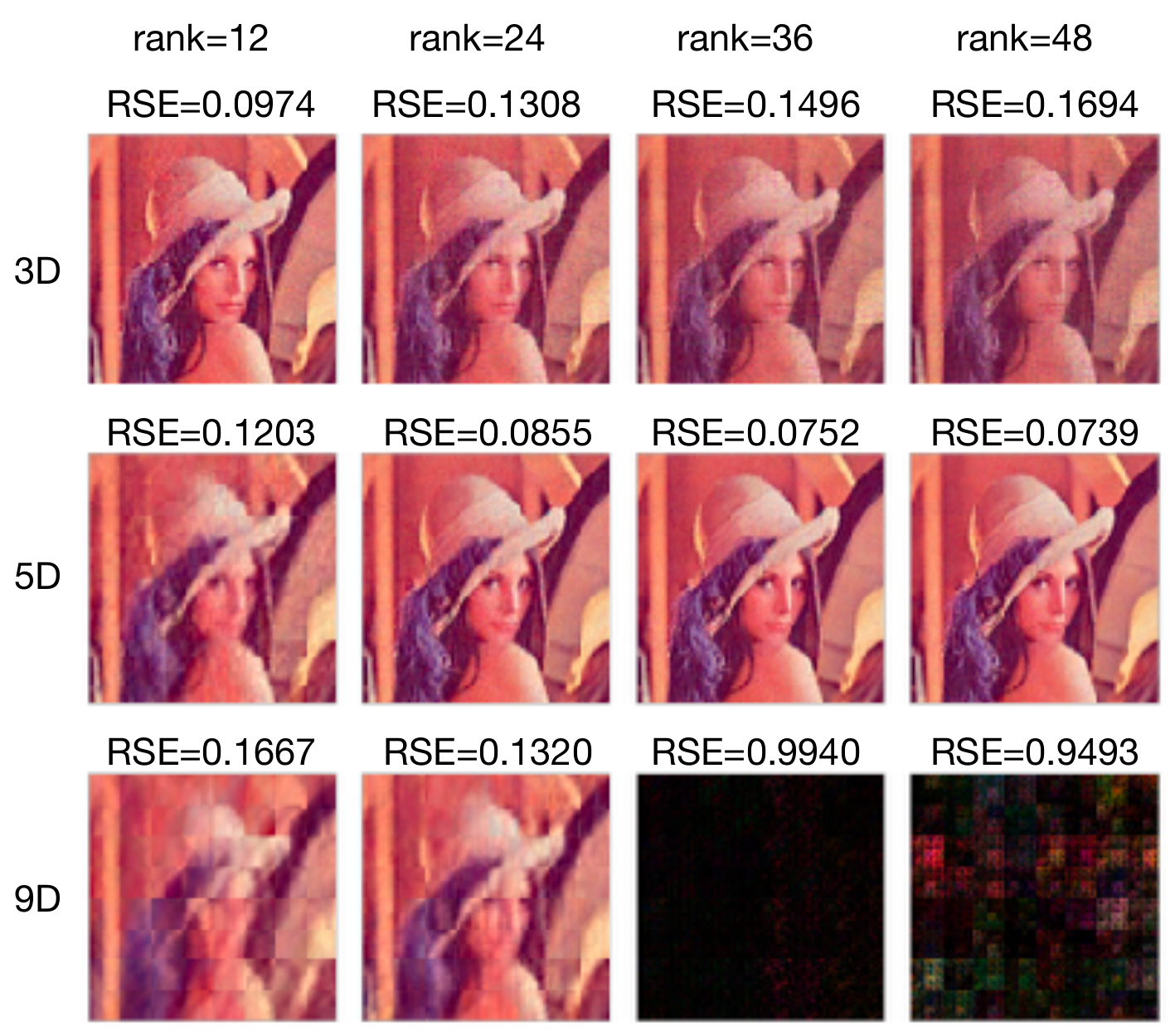}
\caption{Completion results of TR-WOPT under different tensor dimensions and different TR-ranks using benchmark image ``Lena" (missing rate is 0.7).}
\label{lena_rd}
\end{figure}

We choose the benchmark RGB image ``Lena" (original size $256 \times 256 \times 3$) in this experiment and tensorize the image by the above tensorization scheme to 3-D ($256 \times 256 \times 3$), 5-D ($16 \times 16\times 16 \times 16 \times 3$) and 9-D ($4 \times 4 \times 4 \times 4 \times 4 \times 4 \times 4 \times 4 \times 3 $). For simplicity, the TR-ranks of each experiment are set as the same value, i.e., $R_1=R_2=\ldots=R_N$. We conduct the experiments when TR-ranks are 12, 24, 36 and 48 respectively. The missing rates of all the data are set as 0.7. Figure \ref{lena_rd} shows the visual results and corresponding RSE values of each situation. From the results we can see, the best completion performance is obtained at the 5-D case when TR-rank is set as 48. So tensorizing data to a higher-dimension properly can enhance the performance of our algorithm. In the following RGB image experiments, we adopt this tensor structure for our algorithm.

\begin{figure}[!htb]
\centering
\includegraphics[width=0.48\textwidth]{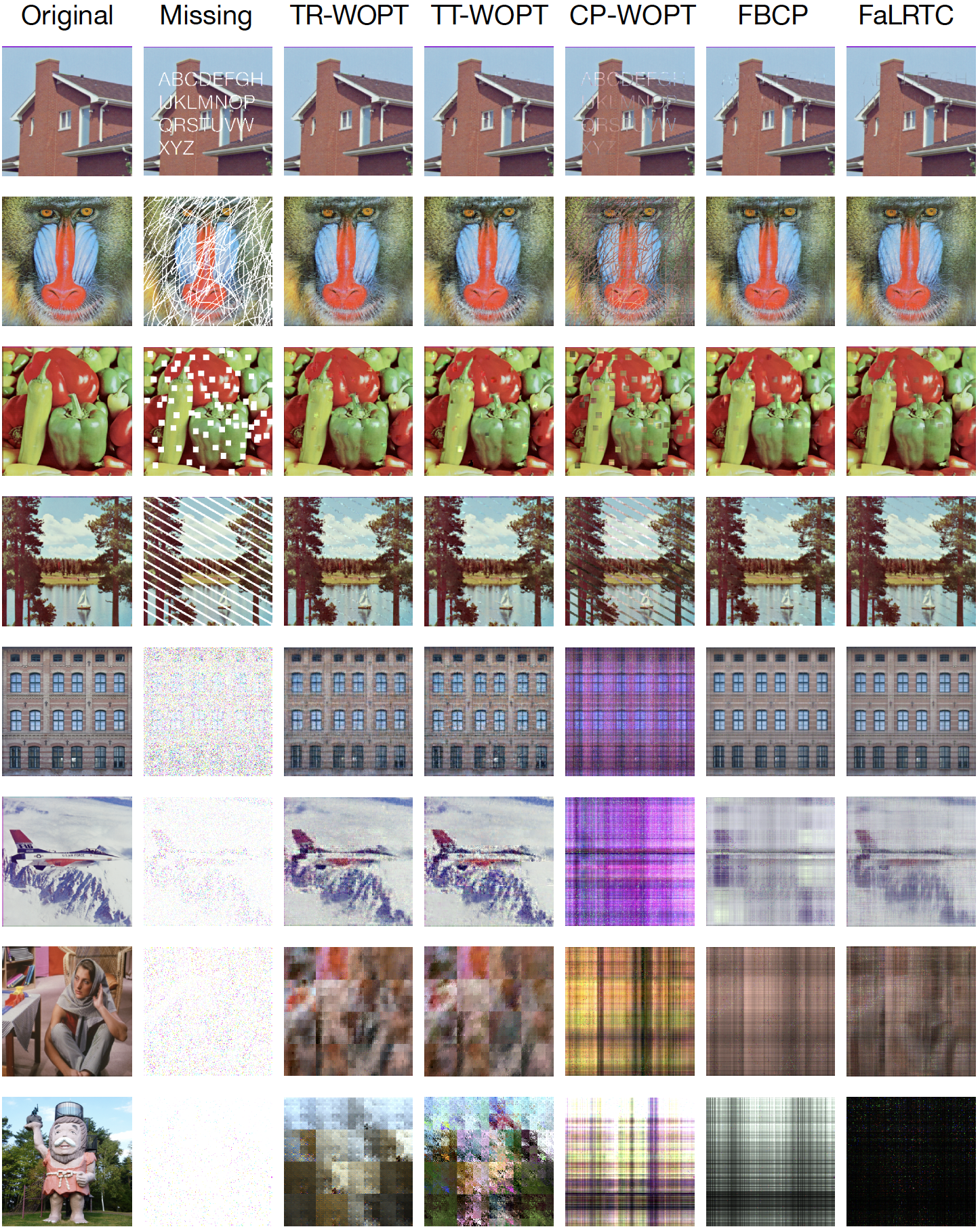}
\caption{Visual completion results of five algorithms under eight image missing situations. The first column and the second column are original images and images with specified missing patterns respectively. The following columns are the completion results of the five algorithms respectively. The first row to the fourth row are the completion results of alphabet missing, scratch missing, block missing and line missing respectively. The fifth row to the last row are random missing completion results of missing rates 0.8, 0.9, 0.95 and 0.99 respectively.}
\label{irr_img}
\end{figure}

The following experiments consider four different irregular missing situations, i.e., removing images by the shapes of the alphabet, missing by scratching, the block missing and the line missing. Moreover, random missing cases with high missing rates are also considered. Because all the tensor completion algorithms perform well in low missing rate situations, we only test high random missing rate situations, i.e., missing rates are 0.8, 0.9, 0.95 and 0.99. We tune ranks and hyper-parameters of each algorithm and record the best completion results of each algorithm. Figure \ref{irr_img} and Table \ref{irr_img_t} show the visual and numerical completion results of the five algorithms respectively. From the results we can see, TR-WOPT performs better than TT-WOPT, CP-WOPT, and FBCP in all the tested situations. However, FaLRTC shows slightly better performance than TR-WOPT in the images of scratch missing and random missing (missing rate is 0.8). This is because the images own distinct low-rank property and the missing rate is relatively low, which is easy for ``rank minimization based" algorithms to catch the low-rank structures of the tensor. However, when the missing rate of data is higher and most of the information is missing, FaLRTC cannot find the low-rank structure of the data, so the performance of random missing rate 0.9 to 0.99 of the algorithm drops quickly.

\begin{table}[!htb]
\centering
\caption{\protect{Numerical completion results of five algorithms under eight image missing situations.}}
\label{irr_img_t}
\setlength{\tabcolsep}{0.8mm}{
\begin{tabular}{c|c|c|c|c|c|c}
\hline
\hline
 \multicolumn{2}{c|}{}& TR-WOPT & TT-WOPT & CP-WOPT & FBCP & FaLRTC\\
 \hline
Alphabet&\makecell[cc]{RSE\\PSNR}   &     \makecell[cc]{ \textbf{0.0227}\\  \textbf{37.17}}    &    \makecell[cc]{ 0.0282\\  35.30}     &
\makecell[cc]{0.0901 \\25.62  }     &           \makecell[cc]{0.0397 \\ 32.32 }   &\makecell[cc]{0.0313 \\  34.40 }\\
\hline
Scratch&\makecell[cc]{RSE \\PSNR}&\makecell[cc]{0.110   \\24.55}&\makecell[cc]{0.119 \\  23.83   }&\makecell[cc]{0.231   \\18.08  }&\makecell[cc]{ 0.114  \\ 24.20 } &   \makecell[cc]{ \textbf{0.106}  \\ \textbf{24.84} }\\
\hline
Block&\makecell[cc]{RSE \\PSNR}&\makecell[cc]{\textbf{0.0891} \\ \textbf{26.21}    }&\makecell[cc]{0.124 \\ 23.31 }&\makecell[cc]{0.176 \\  20.32   }&\makecell[cc]{ 0.115\\ 24.01  }  & \makecell[cc]{ 0.104  \\   24.84  }\\
\hline
Line&\makecell[cc]{RSE \\PSNR}&\makecell[cc]{\textbf{0.101} \\   \textbf{24.81} }&\makecell[cc]{ 0.115  \\  23.70 }&\makecell[cc]{0.187 \\ 19.46  }&\makecell[cc]{0.116 \\ 23.61 }  & \makecell[cc]{0.112 \\ 24.72}\\
\hline
0.8&\makecell[cc]{RSE \\PSNR}&\makecell[cc]{ 0.128   \\ 23.59   }&\makecell[cc]{ 0.142   \\   22.71}&\makecell[cc]{ 0.332  \\   15.32}&\makecell[cc]{ 0.101 \\  25.70}  & \makecell[cc]{\textbf{0.0839}   \\   \textbf{27.27}}\\
\hline
0.9 &\makecell[cc]{RSE \\PSNR}&\makecell[cc]{\textbf{0.125}   \\   \textbf{19.97}}&\makecell[cc]{0.134  \\  19.35}&\makecell[cc]{ 0.414 \\  9.562}&\makecell[cc]{0.175   \\   17.01}  & \makecell[cc]{ 0.146  \\  18.62}\\
\hline
0.95&\makecell[cc]{RSE \\PSNR}&\makecell[cc]{ \textbf{0.252}  \\  \textbf{18.20}}&\makecell[cc]{0.276  \\  17.42}&\makecell[cc]{ 0.530  \\  11.74}&\makecell[cc]{ 0.351 \\  15.32}  & \makecell[cc]{0.343  \\  15.52}\\
\hline
0.99&\makecell[cc]{RSE \\PSNR}&\makecell[cc]{\textbf{0.398}  \\  \textbf{12.59}}&\makecell[cc]{0.657  \\  8.226}&\makecell[cc]{0.679  \\  7.948}&\makecell[cc]{0.457  \\  11.38}  & \makecell[cc]{0.942  \\  5.099}\\
\hline\hline
\end{tabular}
}
\label{irr_num}
\end{table}

\section{Conclusions}
Based on low-rank tensor ring decomposition, in this paper, we proposed a new tensor completion algorithm named tensor-ring weighted optimization (TR-WOPT). The TR core tensors are optimized by the gradient-based method and used to predict the missing entries of the incomplete tensors. We conduct various synthetic data experiments and real-world data experiments, and the results show that TR-WOPT outperforms the state-of-the-art algorithms in many situations. In addition, we also find that tensorizing lower-dimension tensor to a proper higher-dimension tensor can give a better data structure and thus can improve the performance of our algorithm. Good performance of TR-WOPT in various completion tasks shows the high representation ability and flexibility of TR decomposition. It is also shown that the gradient-based algorithm is promising to optimize tensor decompositions. Furthermore, our method needs the TR-rank to be specified before the experiment, which is time-consuming to find the best TR-rank for the data. In our future work, we will study how to determine TR-ranks automatically.
\section*{Acknowledgment}
This work was supported by JSPS KAKENHI (Grant No. 17K00326, 15H04002, 18K04178), JST CREST (Grant No. JPMJCR1784) and the National Natural Science Foundation of China (Grant No. 61773129).

\bibliographystyle{plain}
\bibliography{paper}
\end{document}